\documentclass[12pt,oneside,reqno]{amsart}

\usepackage[a4paper]{geometry}

\usepackage[colorlinks,linkcolor=black,citecolor=blue]{hyperref}
\usepackage{amssymb}
\usepackage[dvips]{graphics}
\usepackage{times}
\usepackage{bbm}
\usepackage{cases}
\usepackage{amsmath}
\usepackage{graphicx}
\usepackage{mathrsfs}
\usepackage{amsfonts}
\usepackage{enumerate,amsmath,amssymb,amsthm}

\numberwithin{equation}{section}

\newcommand{\be}{\begin{eqnarray}}
\newcommand{\ee}{\end{eqnarray}}
\newcommand{\ce}{\begin{eqnarray*}}
\newcommand{\de}{\end{eqnarray*}}
\newtheorem{theorem}{Theorem}[section]
\newtheorem{lemma}[theorem]{Lemma}
\newtheorem{definition}[theorem]{Definition}
\newtheorem{proposition}[theorem]{Proposition}
\newtheorem{Examples}{Example}
\newtheorem{corollary}[theorem]{Corollary}
\theoremstyle{remark}
\newtheorem{remark}[theorem]{Remark}

\def\eps{\varepsilon}

\def\e{\mathrm{e}}

\def\p{\partial}

\def\[{{\Big[}}
\def\]{{\Big]}}
\def\<{{\langle}}
\def\>{{\rangle}}
\def\({{\Big(}}
\def\){{\Big)}}

\def\bx{{\mathbf{x}}}

\def\dif{{\mathord{{\rm d}}}}

\def\={&\!\!=\!\!&}

\def\sL{{\mathcal L}}

\def\mE{{\mathbb E}}

\def\mN{{\mathbb N}}

\def\mP{{\mathbb P}}

\def\mR{{\mathbb R}}

\def\1{{\mathbf{1}}}

\def\sL{{\mathscr L}}

\def\geq{\geqslant}
\def\leq{\leqslant}
\def\ge{\geqslant}

\def\eps{\varepsilon}

\def\e{\mathrm{e}}

\def\p{\partial}

\def\[{{\Big[}}
\def\]{{\Big]}}
\def\<{{\langle}}
\def\>{{\rangle}}
\def\({{\Big(}}
\def\){{\Big)}}

\def\bx{{\mathbf{x}}}

\def\dif{{\mathord{{\rm d}}}}

\def\={&\!\!=\!\!&}
\def\bt{\begin{theorem}}
\def\et{\end{theorem}}
\def\bl{\begin{lemma}}
\def\el{\end{lemma}}
\def\br{\begin{remark}}
\def\er{\end{remark}}
\def\bx{\begin{Examples}}
\def\ex{\end{Examples}}
\def\bd{\begin{definition}}
\def\ed{\end{definition}}
\def\bp{\begin{proposition}}
\def\ep{\end{proposition}}
\def\bc{\begin{corollary}}
\def\ec{\end{corollary}}
\def\bpf{\begin{proof}}
\def\epf{\end{proof}}

\def\geq{\geqslant}
\def\leq{\leqslant}
\def\ge{\geqslant}

\def\<{\langle} 
\def\>{\rangle}

\def\e{\text{\rm{e}}}

\allowdisplaybreaks

\begin{document}

\title[Gradient estimates of heat kernels]{Optimal gradient estimates of heat kernels of stable-like operators}

\author{Kai Du}
\address{Shanghai Center for Mathematical Sciences, Fudan University, Shanghai 200433, P.R.China}
\email{kdu@fudan.edu.cn}

\author{Xicheng Zhang}
\address{School of Mathematics and Statistics, Wuhan University,
Wuhan, Hubei 430072, P.R.China}
\email{XichengZhang@gmail.com}
\thanks{Research of X. Zhang is partially supported by NNSFC grant of China (No. 11731009) and the DFG through the CRC 1283 
``Taming uncertainty and profiting from randomness and low regularity in analysis, stochastics and their applications''. }

\subjclass[2010]{60G52, 35K08.}
\keywords{Gradient estimate, heat kernel, stable-like operators}

\begin{abstract}
In this note we show the optimal gradient estimate for heat kernels of stable-like operators by providing a counterexample.
\end{abstract}

\maketitle

\section{Introduction}

Let $\kappa:\mR^d\to[0,\infty)$ be a bounded measurable function. For $\alpha\in(0,2)$, consider the following nonlocal 
$\alpha$-stable-like operator
$$
\sL_\kappa f(x):=\int_{\mR^d}(f(x+y)-f(x)-y^{(\alpha)}\cdot\nabla f(x))\frac{\kappa(y)}{|y|^{d+\alpha}}\dif y,
$$
where
$$
y^{(\alpha)}:=1_{\alpha\in(1,2)}y+1_{\alpha=1}1_{|y|\leq 1}y.
$$
It is well known that if for some $K_0\geq 1$,
\begin{align}\label{Ka}
K^{-1}_0\leq \kappa(y)\leq K_0,\ \ 1_{\alpha=1}\int_{r<|y|<R}\kappa(y)\dif y=0,\ 0<r<R<\infty,
\end{align}
then there is a smooth fundamental solution $p_\kappa(t,x)$ to the operator $\sL_\kappa$ satisfying (see \cite{Ch-Zh, Ch-Zh2, Ji1})
$$
\p_tp_\kappa(t,x)=\sL_\kappa p_\kappa(t,\cdot)(x),\ t>0,x\in\mR^d.
$$ 
Moreover, $p_\kappa(t,x)$ enjoys the following two-sided estimates: for some $K_1=K_1(\alpha,d,K_0)\geq 1$,
$$
K_1^{-1}t(t^{1/\alpha}+|x|)^{-d-\alpha}\leq p_\kappa(t,x)\leq K_1t(t^{1/\alpha}+|x|)^{-d-\alpha},
$$
and gradient estimate: for some $K_2=K_2(\alpha,d,K_0)\geq 1$,
\begin{align}\label{grad}
|\nabla p_\kappa(t,x)|\leq K_2 t^{1-1/\alpha}(t^{1/\alpha}+|x|)^{-d-\alpha}.
\end{align}
The above estimates can be found in \cite{Ch-Zh2, Ji1}.
Notice that for $\lambda>0$, if we let $\kappa_\lambda(y):=\kappa(\lambda^{1/\alpha}y)$, then $p_\kappa$ has the following scaling property:
$$
p_\kappa(\lambda t,\lambda^{1/\alpha}x)=\lambda^{-d/\alpha}p_{\kappa_\lambda}(t,x).
$$
Moreover, when $\kappa(y)=1$, it is well known that
$$
\sL_1=c_{\alpha,d}\,\Delta^{\alpha/2},
$$
where $c_{\alpha,d}>0$ and $\Delta^{\alpha/2}:=-(-\Delta)^{\alpha/2}$ is the usual fractional Laplacian.
In this case, it is also well known that the sharp gradient estimate takes the following form: for some $K_3=K_3(\alpha,d)\geq 1$,
\begin{align}\label{Sharp}
|\nabla_{\!x} p_1(t,x)|\leq K_3 t(t^{1/\alpha}+|x|)^{-d-\alpha-1}.
\end{align}
Indeed, let $S_t$ be the $\alpha/2$-subordinator and $\phi(t,x)=(2\pi t)^{-d/2}\e^{-|x|^2/(2t)}$ be the Gaussian heat kernel. 
By the subordination we have
\begin{align}\label{He}
p_1(t,x)=\int^\infty_0\phi(s,x)\,\mP\circ S^{-1}_t(\dif s).
\end{align}
Since $\mP\circ S^{-1}_t(\dif s)\leq C ts^{-1-\alpha/2}\e^{-ts^{-\alpha/2}}\dif s$, by elementary calculations, it follows that (see \cite{Bo-Ja})
$$
|\nabla_{\!x} p_1(t,x)|\leq Ct|x|\int^\infty_0s^{-2-(d+\alpha)/2}\e^{-ts^{-\alpha/2}-|x|^2/(2s)}\dif s\leq K_3 t(t^{1/\alpha}+|x|)^{-d-\alpha-1}.
$$
Notice that the right hand side of \eqref{Sharp} is smaller than the one in \eqref{grad} when $x$ goes to $\infty$. 
We mention that gradient estimate \eqref{Sharp} plays an important role 
in \cite[Proposition 3.2]{Kn-Ku}.

Here a natural question is that for general stable-like operator $\sL_\kappa$, is it possible to show the same gradient estimate \eqref{Sharp}?
We have the following negative answer. Thus \eqref{grad} is optimal.
\bt\label{Main}
For $d=1$, there is an even function $\kappa:\mR\to[1,2]$ such that
$$
\limsup_{|x|\to\infty}|x|^{1+\alpha}|\nabla_{\!x} p_\kappa(1,x)|>0.
$$
\et
\br
Although we construct a symmetric example, it is in fact also easier to construct  non-symmetric examples from the following proofs.
\er

\section{Proof of Theorem \ref{Main}}
Below we assume $d=1$, and simply write
$$
p_\kappa(x)=p_\kappa(1,x).
$$
Suppose that $\kappa$ satisfies \eqref{Ka} and $f:\mR\to[0,1]$ is a measurable function so that
$$
\lambda:=\int_\mR f(x)|x|^{-1-\alpha}\dif x<\infty.
$$
Let $N^\lambda_t$ be a Poisson process with intensity $\lambda$ and
$\xi_1,\cdots\xi_n,\cdots$ a sequence of i.i.d. random variables with common distribution density $q(x):=f(x) |x|^{-1-\alpha}/\lambda$, 
which are independent with $N^\lambda$.
Let $\xi_0:=0$ and
$$
Y^f_t:=\xi_0+\xi_1+\cdots+\xi_{N^\lambda_t}.
$$
Then $Y^f_t$ is a compound Poisson process with L\'evy exponent
$$
\int_{\mR}(\e^{i\xi\cdot y}-1) f(y)|y|^{-1-\alpha}\dif y.
$$
Let $Z^\kappa_t$ be an independent L\'evy process with L\'evy measure $\kappa(y)|y|^{-1-\alpha}\dif y$. The sum $Z^\kappa_t+Y^f_t$
is still a L\'evy process with L\'evy measure $(\kappa(y)+f(y))|y|^{-1-\alpha}\dif y$. In particular, we have
\begin{align}\label{Es5}
p_{\kappa+f}(x)=\mE p_\kappa(x-Y^f_1).
\end{align}
\bl\label{Le21}
For $a,\eps>0$, let $f^\eps_a(x):=f\big(\frac{a+x}{\eps}\big)+f\big(\frac{a-x}{\eps}\big)$. 
Suppose that 
$$
{\rm supp}(f)\subset[-1,1],\ \ \int_\mR f(x)\,\dif x=1,
$$
and for some $z_0\in(-1,0)$, $\eps,\delta\in(0,1)$ and $\gamma, A\geq 1$,
\begin{align}\label{Con1}
\inf_{x\in[z_0-\eps,z_0+\eps]}p'_\kappa(x)\geq\delta,\quad\|p'_\kappa\|_\infty\leq\gamma,\quad |p'_\kappa(x)|\leq \gamma|x|^{-2-\alpha} \ \text{ for }\ x\geq A.
\end{align}
Then there are constants $C=C(\alpha,\delta,\gamma)>0$ and $A_0=A_0(\delta,\gamma)\geq 2$ such that for all $a\geq A\vee A_0$,
\begin{align}\label{Es6}
p'_{\kappa+f^\eps_a}(z_0+a)\geq C a^{-1-\alpha}.
\end{align}
\el
\begin{proof}
For $a>1$, let 
$$
\lambda:=\int_\mR f^\eps_a(x)|x|^{-1-\alpha}\dif x=2\int_{\mR}f\Big(\frac{a+x}{\eps}\Big)|x|^{-1-\alpha}\dif x=
2\int_{\mR}f\Big(\frac{a-x}{\eps}\Big)|x|^{-1-\alpha}\dif x.
$$
By \eqref{Es5} and the definition of $Y_1^{f^\eps_a}$, we have
\begin{align*}
p'_{\kappa+f^\eps_a}(z_0+a)&=p'_\kappa(z_0+a)\mP(N^\lambda_1=0)+\mE p'_\kappa(z_0+a-\xi_1)\mP(N^\lambda_1=1)\\
&\quad +\mE p'_\kappa(z_0+a-Y^{f^\eps_a}_1+\xi_1)\mP(N^\lambda_1\geq 2)\\
&=:J_0+J_1+J_2.
\end{align*}
For $J_0$, by \eqref{Con1} we have
$$
|J_0|\leq |p'_\kappa(z_0+a)|\e^{-\lambda}\leq \gamma|z_0+a|^{-2-\alpha},\ a\geq A+1.
$$
For $J_1$, since ${\rm supp}(f)\subset[-1,1]$, we have
\begin{align*}
J_1&=\e^{-\lambda}\int_{\mR}p'_\kappa(z_0+a+y)\Big(f\Big(\frac{a+y}{\eps}\Big)+f\Big(\frac{a-y}{\eps}\Big)\Big) |y|^{-1-\alpha}\dif y\\
&\geq \e^{-\lambda}\int_{\mR}\left(\delta\cdot f\Big(\frac{a+y}{\eps}\Big)-\gamma(2a-\eps-|z_0|)^{-2-\alpha}f\Big(\frac{a-y}{\eps}\Big)\right)|y|^{-1-\alpha}\dif y\\
&=\e^{-\lambda}\big(\delta-\gamma(2a-\eps-|z_0|)^{-2-\alpha}\big)\lambda/2.
\end{align*}
For $J_2$, we have
\begin{align*}
|J_2|\leq \|p'_\kappa\|_\infty(1-\e^{-\lambda}-\lambda\e^{-\lambda})\leq 2\gamma\lambda^2,\ \lambda\in(0,1).
\end{align*}
Combining the above calculations, and thanks to
$$
2\eps(a+\eps)^{-1-\alpha}\leq \lambda\leq 2\eps(a-\eps)^{-1-\alpha},
$$
we obtain that for $a\geq (2\gamma/\delta)^{2+\alpha}\vee 2$,
\begin{align*}
p'_{\kappa+f^\eps_a}(z_0+a)\geq \lambda\e^{-\lambda}\delta/4-\gamma|z_0+a|^{-2-\alpha}-2\gamma\lambda^2\geq C a^{-1-\alpha}.
\end{align*} 
Thus we complete the proof.
\end{proof}

Recalling that $p_1(x)=p_1(1,x)$ is given by \eqref{He}, in the following we may fix $z_0\in(-1,0)$ and $\eps\in(0,|z_0|/2)$ so that
\begin{align}\label{Def}
\delta:=\inf_{x\in[z_0-\eps,z_0+\eps]}p'_1(x)\in(0,1).
\end{align}
In particular, by Lemma \ref{Le21} we have
\bc
There are constants $C=C(\alpha,\delta, z_0)>0$ and $A>0$ such that for any $a>A$, one can find an even function $\kappa:\mR\to[1,2]$ 
depending on $a$ such that
\begin{align}\label{Es6}
p'_\kappa(z_0+a)\geq Ca^{-1-\alpha}.
\end{align}
\ec

The above corollary implies that it is not possible to find a constant $C$ that only depends on the bound of $\kappa$ so that for all 
$\kappa:\mR\to[1,2]$ and $x\in\mR$,
$$
|p'_\kappa(x)|\leq C (1+|x|)^{-2-\alpha}.
$$
In fact, it has already given a negative answer to our question.
Nevertheless, we are still interested in finding an even function $\kappa:\mR\to[1,2]$ so that \eqref{Es6} holds for some sequence $a_n\to\infty$.
We prepare the following lemma.

\bl\label{Le23}
Let $\kappa(x)=1+f(x)$, where $f:\mR\to[0,1]$ is an even function so that
$$
\lambda:=\int_\mR f(x)|x|^{-1-\alpha}\dif x<\infty.
$$
(i) Under \eqref{Def}, there is a $\lambda_0=\lambda_0(\alpha,\delta)>0$ such that for all the above $f$ with $\lambda\leq \lambda_0$,
\begin{align}\label{Es1}
\inf_{x\in[z_0-\eps,z_0+\eps]}p'_\kappa(x)\geq\delta/2.
\end{align}
(ii) There exists a constant $\gamma=\gamma(\alpha,d)>0$ such that for all $A\geq 1$ and all the above $f$ with ${\rm supp}(f)\subset[-A,A]$ and $\lambda\leq 1$,
\begin{align}\label{Es2}
|p'_\kappa(x)|\leq \gamma|x|^{-2-\alpha},\ \ |x|\geq A^2.
\end{align}
\el
\begin{proof}
(i) By \eqref{Es5} with $\kappa=1$ there, we have
\begin{align*}
p'_{1+f}(x)&=\mP(N^\lambda_1=0) p_1'(x)+\sum_{k\geq 1}\mP(N^\lambda_1=k)\int_\mR p'_1(x-y)q^{*k}(y)\dif y\\
&\geq \e^{-\lambda} p'_1(x)-\|p'_1\|_\infty\sum_{k\geq 1}\mP(N^\lambda_1=k)\\
&=\e^{-\lambda} \Big(p'_1(x)-\|p'_1\|_\infty(\e^{\lambda}-1)\Big),
\end{align*}
where $q(x):=f(x)|x|^{-1-\alpha}/\lambda$ and $q^{*k}$ stands for the $k$-order convolution.
In particular, one can choose $\lambda_0$ small enough so that
$$
(1\vee\|p'_1\|_\infty)(\e^{\lambda_0}-1)\leq\delta/3.
$$
Estimate \eqref{Es1} then follows by definition \eqref{Def}.
\\
\\
(ii) We make the following decomposition:
\begin{align*}
p'_\kappa(x)&=\Bigg(\sum_{k\leq \sqrt{|x|}/2}+\sum_{k>\sqrt{|x|}/2}\Bigg)\mP(N^\lambda_1=k)\int_\mR p'_1(x-y)q^{*k}(y)\dif y=:J_1(x)+J_2(x).
\end{align*}
For $J_1$, noticing that ${\rm supp}(q^{*k})\subset[-kA,kA]$, by \eqref{Sharp} and $|x|>A^2$, we have
\begin{align*}
|J_1(x)|&\leq \sum_{k\leq \sqrt{|x|}/2}\mP(N^\lambda_1=k)\int_\mR |p'_1(x-y)| q^{*k}(y)\dif y\\
&\leq K_3\sum_{k\leq \sqrt{|x|}/2}\mP(N^\lambda_1=k)\Big(|x|-\sqrt{|x|}A/2\Big)^{-2-\alpha}\\
&\leq K_3 2^{2+\alpha}|x|^{-2-\alpha}.
\end{align*}
For $J_2$, by Stirling's formula and $\lambda\leq 1$, we have
$$
|J_2(x)|\leq \|p'_1\|_\infty\e^{-\lambda}\sum_{k>\sqrt{|x|}/2}\lambda^k/k!\leq C|x|^{-2-\alpha}.
$$
Combining the above calculations, we obtain \eqref{Es2}.
\end{proof}

Now we are in a position to give
\begin{proof}[Proof of Theorem \ref{Main}]
Let $h:\mR\to[0,1]$ be a measurable function with
$$
{\rm supp}(h)\subset[-1,1],\quad \int_\mR h(x)\dif x=1.
$$
Let $A\geq A_0$ be a large number, whose value will be determined later, where $A_0$ is from Lemma \ref{Le21}. For $k\in\mN$, define
\begin{align*}
 A_1 & := A,\quad A_{k+1}:= (A_k + \eps)^2,\ k \ge 1,\\
h_k(x) & :=h\Big(\frac{A_k-x}{\eps}\Big)+h\Big(\frac{A_k+x}{\eps}\Big),\quad  \kappa(x):=1+\sum_{k=1}^\infty h_k(x).
\end{align*}
Clearly, $\kappa(-x)=\kappa(x)$. We want to show that for some $C_0=C_0(\alpha)>0$,
\begin{align}\label{Es8}
p'_{\kappa}(z_0+A_{n})\geq C_0A_{n}^{-1-\alpha},\ \forall n\in\mN.
\end{align}
First of all, by definition we have
\begin{align}\label{Es7}
{\rm supp}(h_k)\subset[A_k-\eps,A_k+\eps]\cup[-A_k-\eps,-A_k+\eps]
\end{align}
and
\begin{align*}
2\eps(A_k+\eps)^{-1-\alpha} 
& \leq \int_\mR h_k(x)|x|^{-1-\alpha}\dif x \nonumber\\
& =2\int_\mR h\Big(\frac{A_k-x}{\eps}\Big)|x|^{-1-\alpha}\dif x\leq 2\eps(A_k-\eps)^{-1-\alpha}.
\end{align*}
For $n\in\mN$, define
$$
\kappa_n(x):=1+\sum_{k=1}^n h_k(x)=:1+f_n(x).
$$
Noticing that ${\rm supp}(f_n)\subset \cup_{k=1}^n {\rm supp}(h_k)\subset[-A_n-\eps, A_n+\eps]$ and
$$
\int_\mR f_n(x)|x|^{-1-\alpha}\dif x\leq 2\eps\sum_{k=1}^n(A_k-\eps)^{-1-\alpha}\leq C A^{-1-\alpha},
$$
by Lemma \ref{Le23}, if $CA^{-1-\alpha}\leq \lambda_0$, then 
\begin{align}\label{Es3}
\inf_{x\in[z_0-\eps,z_0+\eps]}p'_{\kappa_n}(x)\geq\delta/2,
\end{align}
and for some $\gamma=\gamma(\alpha)\geq 1$,
\begin{align}\label{Es4}
|p'_{\kappa_n}(x)|\leq \gamma|x|^{-2-\alpha},\ \ |x|  \ge (A_n+\eps)^2 = A_{n+1}.
\end{align}
Using these two estimates, and by Lemma \ref{Le21} we derive that for some $C_1=C_1(\alpha,\delta)>0$,
$$
p'_{\kappa_{n+1}}(z_0+A_{n+1})=p'_{\kappa_{n}+h_{n+1}}(z_0+A_{n+1})\geq C_1A_{n+1}^{-1-\alpha}.
$$
Finally, let $\tilde f_n(x):=\kappa(x)-\kappa_n(x)=\sum_{k>n}h_k(x)$ and
$$
\tilde \lambda:=\int_\mR\tilde f_n(x)|x|^{-1-\alpha}\dif x\leq 2\eps\sum_{k=n+1}^\infty(A_k-\eps)^{-1-\alpha}\leq C_2A^{-1-\alpha}_{n+1}. 
$$
By \eqref{grad}, there is a constant $M>0$ such that for all $n\in\mN$,
$$
|p'_{\kappa_n}(x)|\leq M,\ \ x\in\mR^d.
$$
As above, we have
\begin{align*}
p'_{\kappa}(z_0+A_{n})&=\mP(N^{\tilde \lambda}_1=0) p'_{\kappa_n}(z_0+A_n)
+\sum_{k\geq 1}\mP(N^{\tilde\lambda}_1=k)\int_\mR p'_{\kappa_n}(z_0+A_n-y){\tilde q}^{*k}(y)\dif y\\
&\geq \e^{-\tilde\lambda} p'_{\kappa_n}(z_0+A_n)-M\sum_{k\geq 1}\mP(N^{\tilde\lambda}_1=k)\\
&=\e^{-\tilde\lambda}(p'_{\kappa_n}(z_0+A_n)-M(\e^{\tilde\lambda}-1))\\
&\geq \e^{-\tilde\lambda}(C_1A_n^{-1-\alpha}-2C_2MA^{-1-\alpha}_{n+1})\geq C_3A_n^{-1-\alpha},
\end{align*}
provided $A$ large enough.
Thus we get \eqref{Es8}, which means 
$$
\limsup_{x\to\infty}x^{1+\alpha}p'_\kappa(x)\geq C_0>0.
$$
By symmetry, we obtain the desired estimate.
\end{proof}

{\bf Acknowledgement:}
The authors would like to thank Peng Jin and Guohuan Zhao for their useful discussions.


\begin{thebibliography}{10}

\bibitem{Bo-Ja} K. Bogdan, T. Jakubowski. 
\newblock Estimates of heat kernel of fractional Laplacian perturbed by gradient operators. {\em Comm.
Math. Phys.} {\bf 271} (1) (2007) 179-198.

 \bibitem{CK08}  Z.-Q. Chen and T. Kumagai,
\newblock Heat kernel estimates for  jump processes of mixed types on metric measure spaces.
{\em  Probab. Theory Related Fields \bf 140} (2008), 277-317.

\bibitem{Ch-Zh}
Z.-Q. Chen and X.~Zhang.
\newblock Heat kernels and analyticity of non-symmetric jump diffusion
  semigroups.
\newblock {\em Probab. Theory Related Fields}, {\bf 165}(1-2):267--312, 2016.

\bibitem{Ch-Zh2}
Z.-Q. Chen and X.~Zhang.
\newblock  Heat kernels for time-dependent non-symmetric stable-like operators.
{\it J. Math. Anal. Appl.} {\bf 465} (2018) 1-21.

\bibitem{Ji1}P. Jin. 
\newblock Heat kernel estimates for non-symmetric stable-like processes. https://arxiv.org/abs/1709.02836.
 
\bibitem{Kn-Ku} V. Knopova and A. Kulik. 
\newblock Parametrix construction of the transition probability density of the solution to an SDE driven by $\alpha$-stable noise.
https://arxiv.org/abs/1412.8732v2.
 \end{thebibliography}
\end{document}